\documentclass{amsart}

\newtheorem{theorem}{Theorem}[section]
\newtheorem{lemma}[theorem]{Lemma}

\theoremstyle{definition}
\newtheorem{definition}[theorem]{Definition}

\theoremstyle{remark}

\numberwithin{equation}{section}

\newcommand{\abs}[1]{\lvert#1\rvert}

\begin{document}

\title{Zeros of Meromorphic function}


\author{Lande Ma}
\curraddr{School of Mathematical Sciences, Tongji University, Shanghai, 200092, China}
\email{dzy200408@126.com}
\thanks{}

\author{Zhaokun Ma}
\address{}
\curraddr{YanZhou College, ShanDong Radio and TV University, YanZhou, ShanDong 272100 China}
\email{dzy200408@sina.cn}
\thanks{}


\date{April 19th, 2022.}

\dedicatory{}


\begin{abstract}
We introduce and develop the root locus method in mathematics. And we study the distribution of zeros of meromorphic functions by root locus method.
\end{abstract}

\maketitle

\section*{Introduction}

The meromorphic functions are the basic functions with general meanings in mathematics. Meromorphic functions include all complex functions in the engineering. The complex functions that need to be studied in the engineering must have poles. The poles mainly determine the dynamic characteristics of a process. Complex functions without poles have no actual values in the engineering. The main purpose of mathematical research is to provide tools for engineering and other sciences. In pure mathematics, many functions like the Riemann zeta function and gamma function, belong to the meromorphic functions. However, there are only several properties concerning the general meromorphic functions. And, there is almost no result of its derivative function. Regarding the properties of zeros of derivative of general meromorphic functions, so far, except for the definition of the zero of derivative, nothing has been discovered.

It is well known that the distribution of zeros of the derivative of meromorphic functions has been studied\cite{Stein,Bergweiler,Yamanoi}, and several special meromorphic functions are also studied due to their great influence, such as Riemann zeta function\cite{Berndt,Levinson,Conrey,Zhang,Ge,Edwards}. By utilizing the results which we obtained, we get the necessary and sufficient results about the zeros of derivatives of meromorphic functions.

Speiser showed that the Riemann hypothesis is equivalent to the absence
of non-trivial zeros of the derivative of the Riemann zeta-function left of the critical line\cite{Speiser}.
Levinson and Montgomery proved the quantitative version of the Speiser's result, namely, that the Riemann zeta-function and
its derivative have approximately the same number of zeros left of the critical line\cite{Levinson}.
We prove there don't exist the zeros of derivatives of Riemann zeta function and Dirichlet ate series on the critical line or on the imaginary axis. We also prove that the results of the values of argument of Riemann zeta function and Dirichlet ate series strict and monotonic decrease on the critical line. Riemann hypothesis is equivalent to that argument of xi-function and Riemann zeta function and Dirichlet ate series concerning the imaginary variable $t$ strictly and monotonically decrease at the left side of the critical line.

\section{The results }
In textbook of \emph{automatic control theory}\cite{Richard}, the factor in the left of the root locus equation is the rational fraction function of the constant coefficient. The root locus equation is only concerning two degree of 0 and 180 degree which obtains the real number values. So, the root locus equations and the results of the root locus in \emph{automatic control theory} are all very special and limit. The proofs of results in \emph{automatic control theory} are not comprehensive and no accurate. For the pure mathematical study, the factor which is in the left side of the root locus equation is best substituted by a general meromorphic function. So, we need to give the results of the root locus of the general meromorphic function.

If the number of poles is finite or infinite, and the number of zeros is finite or infinite, namely the two finite or infinite products $\prod^{m}_{l=1}(1-\frac{s}{z_{l}})^{\gamma_{l}}G_{lz}(s)$ and $\prod^{n}_{j=1}(1-\frac{s}{p_{j}})^{\beta_{j}}G_{jp}(s)$ both converge, and $K$ is as: $K=|\frac{\prod^{n}_{j=1}(1-\frac{s}{p_{j}})^{\beta_{j}}G_{jp}(s)}
{G(s)\prod^{m}_{l=1}(1-\frac{s}{z_{l}})^{\gamma_{l}}G_{lz}(s)}|$. The $K$ is the reciprocal of the modulus of the meromorphic function $G(s)\prod^{m}_{l=1}(1-\frac{s}{z_{l}})^{\gamma_{l}}G_{lz}(s)/\prod^{n}_{j=1}(1-\frac{s}{p_{j}})^{\beta_{j}}G_{jp}(s)$, let, $W(s)=\frac{G(s)\prod^{m}_{l=1}(1-\frac{s}{z_{l}})^{\gamma_{l}}G_{lz}(s)}
{\prod^{n}_{j=1}(1-\frac{s}{p_{j}})^{\beta_{j}}G_{jp}(s)}=u(\sigma,t)+iv(\sigma,t)$, $W(s)$ is a meromorphic function with one complex variable in the extended complex plane $\mathbb{C}\cup\{\infty\}$. So, after $K$ and the  meromorphic function $W(s)$ are multiplied, the product result is the unit complex value of the meromorphic function.

\begin{definition}
In the extended complex plane $\mathbb{C}\cup\{\infty\}$, the equation (1.1) is called as \emph{the root locus equation}.
\begin{equation}
KG(s)\frac{\prod^{m}_{l=1}(1-\frac{s}{z_{l}})^{\gamma_{l}}G_{lz}(s)}
{\prod^{n}_{j=1}(1-\frac{s}{p_{j}})^{\beta_{j}}G_{jp}(s)}=a+ib
\end{equation}
\end{definition}
In which, $a+ib$ is the unit complex number value of meromorphic function $W(s)$. The argument of the unit complex number value is as: $\alpha=2q\pi+\arg(\frac{b}{a})$. The factor at the left side of the root locus equation (1.1) is meromorphic function $W(s)$.

The zeros $z_{l}$ and poles $p_{j}$ are points in $\mathbb{C}\cup\{\infty\}$, and may be no conjugate. The exponents $\gamma_{l}$ and $\beta_{j}$ are all non-negative real numbers, and they are the exponents of zeros $z_{l}$ and poles $p_{j}$ of Eq(1.1).
\begin{lemma}
After the coincident and finite zeros and poles of Eq(1.1) are cancelled, for any finite zero of Eq(1.1), its $K$ value is $K=+\infty$; for any finite pole of Eq(1.1), its $K$ value is $K=0$.
\end{lemma}
\begin{proof}
The Eq(1.1) can be transformed to the following characteristic equation.

$KG(s)\prod^{m}_{l=1}(1-\frac{s}{z_{l}})^{\gamma_{l}}G_{lz}(s)-
(a+ib)\prod^{n}_{j=1}(1-\frac{s}{p_{j}})^{\beta_{j}}G_{jp}(s)=0$. Obviously, if $K=0$, all roots of the pole factor $\prod^{n}_{j=1}(1-\frac{s}{p_{j}})^{\beta_{j}}G_{jp}(s)$ are the roots of characteristic equation. Conversely, all roots of characteristic equation if $K=0$ are all  roots of the pole factor. So, at poles of Eq(1.1), $K=0$.

The Eq(1.1) can be transformed to its another characteristic equation.
$G(s)\prod^{m}_{l=1}(1-\frac{s}{z_{l}})^{\gamma_{l}}G_{lz}(s)-
\frac{a+ib}{K}\prod^{n}_{j=1}(1-\frac{s}{p_{j}})^{\beta_{j}}G_{jp}(s)=0$. If $K=+\infty$, all roots of the zero factor $\prod^{m}_{l=1}(1-\frac{s}{z_{l}})^{\gamma_{l}}G_{lz}(s)$ are the roots of the characteristic equation. Conversely, all roots of the characteristic equation when $K=+\infty$ are the roots of the zero factor. So, at zeros of Eq(1.1), $K=+\infty$.
\end{proof}

According to the expression of $K$ of the root locus equation (1.1) and according to Lemma 1.2, it is obvious that the $K$ values are continuous concerning the complex variable $s$ in $\mathbb{C}\cup\{\infty\}$.
So, we can prove Theorem 1.3 simply.

\begin{theorem}
Let $s\in\mathbb{C}\cup\{\infty\}$, the $K$ value of $s\in \mathbb{C}\cup\{\infty\}$ takes all of the non-negative real number value from 0 to the $+\infty$.
\end{theorem}
For the argument\cite{Stein} in complex analysis, in \emph{automatic control theory}, it is called as the phase angle. $K$ is called as \emph{the gain $K$}.

Let $\Delta=\{s=\sigma+it\in\mathbb{C}\cup\{\infty\}$, $s$ is not zero or pole of $W(s)$ and $Eq(1.1)\}$.
\begin{lemma}
For any point $s\in\Delta$, the phase angle of meromorphic function $W(s)$ can be written as the $\varphi(\sigma,t)$ function:
$\varphi(\sigma,t)=\arg(\frac{G_{y}(\sigma,t)}{G_{x}(\sigma,t)})+
\sum^{m}_{l=1}(\gamma_{l}\arg(\frac{t-t_{l}}{\sigma-\sigma_{l}})-\gamma_{l}\arg(\frac{t_{l}}{\sigma_{l}})+
\arg(\frac{G_{zy}(\sigma,t)}{G_{zx}(\sigma,t)}))
-\sum^{n}_{j=1}(\beta_{j}\arg(\frac{t-t_{j}}{\sigma-\sigma_{j}})-\beta_{j}\arg(\frac{t_{j}}{\sigma_{j}})+
\arg(\frac{G_{py}(\sigma,t)}{G_{px}(\sigma,t)}))$.
\end{lemma}
Let $s=(\sigma,t)\in\Delta$, then the phase angle condition equation of Eq(1.1) on $s$ is:
\begin{equation}
\varphi(\sigma,t)=2q\pi+arg{(b/a)}
\end{equation}
In which, $q$ is an integer number.
\begin{lemma}
Let Eq(1.1) be the root locus equation of $2q\pi+\alpha$ degree, for any point $s=(\sigma,t)\in\Delta$, if the point $s$ satisfies the phase angle condition equation (1.2) of $2q\pi+\alpha$ degree, then it must be on the root locus of $2q\pi+\alpha$ degree of Eq(1.1).
\end{lemma}
\begin{proof}
Assuming that a point $s_1=(\sigma_1,t_1)\in\Delta$ satisfies the phase angle condition
equation (1.2), $\varphi(\sigma_1,t_1)=2q\pi+arg{(b_1/a_1)}$. The phase angle of $W(s)$ is $\varphi(\sigma,t)$, this phase
angle expression is same as the expression of the left side of
Eq(1.2).

Because the point $s_{1}$ satisfies Eq(1.2), when the point $s_{1}$ is substituted into the
function $W(s)$, $W(s)$ obtains a complex value, according to the
phase angle expression of Eq(1.2) of the point $s_{1}$ and the
condition which a point $s_{1}$ satisfies the phase angle condition
equation (1.2), this complex value can be written
as: $K_1^{*}(a_{1}+ib_{1})$, and  $\alpha_1=\arg{(b_1/a_1)}$,
$K_1^{*}$ is the modulus of function $W(s)$ of point $s_{1}$.
$K_1^{*}$ is a non-zero positive real value. $\alpha_1$ is the phase
angle of the function $W(s)$ of the point $s_{1}$, and it is the value
of the right side of Eq(1.2).

If we bring the point $s_1$ into the gain expression $|\frac{\prod^{n}_{j=1}(1-\frac{s}{p_{j}})^{\beta_{j}}G_{jp}(s)}
{G(s)\prod^{m}_{l=1}(1-\frac{s}{z_{l}})^{\gamma_{l}}G_{lz}(s)}|$,  we can obtain
an unique gain $K_1$, $K_1=|\frac{\prod^{n}_{j=1}(1-\frac{s_1}{p_{j}})^{\beta_{j}}G_{jp}(s_1)}
{G(s_1)\prod^{m}_{l=1}(1-\frac{s_1}{z_{l}})^{\gamma_{l}}G_{lz}(s_1)}|$, the gain $K_1$ is a
reciprocal of the modulus $K_1^{*}$ of $W(s)$ of the point $s_1$.
Further obtain $K_1^{*}$. $K_1^{*}=1/K_1$. For that the gain $K_1$
multiplied the expression $G(s_1)\prod^{m}_{l=1}(1-\frac{s_1}{z_{l}})^{\gamma_{l}}G_{lz}(s_1)/\prod^{n}_{j=1}(1-\frac{s_1}{p_{j}})^{\beta_{j}}G_{jp}(s_1)$, we have,
$K_1G(s_1)\prod^{m}_{l=1}(1-\frac{s_1}{z_{l}})^{\gamma_{l}}G_{lz}(s_1)/\prod^{n}_{j=1}(1-\frac{s_1}{p_{j}})^{\beta_{j}}G_{jp}(s_1)=a_1+ib_1$. This equation is a concrete situation which a point $s_1$ satisfies Eq(1.1). The previous results prove
that the point $s_1$ is a root of Eq(1.1).

So, the point $s_1$ satisfies Eq(1.1). Therefore, it is proved
that the point $s_1$  is on the root locus of Eq(1.1), and which its gain is $K_1$ and its degree is
$2q\pi+\alpha_1$. Hence, if the point $s_1$ satisfies Eq(1.2), then
the point $s_1$ is on the root locus of Eq(1.1), and which its
gain is $K_1$ and its degree is $2q\pi+\alpha_1$.
\end{proof}
\begin{lemma}
Let Eq(1.1) be the root locus equation of $2q\pi+\alpha$ degree, for any point $s=(\sigma,t)\in\Delta$, if the point $s$ is on the root locus of $2q\pi+\alpha$  degree of Eq(1.1), then it must satisfy the phase angle condition equation
of $2q\pi+\alpha$ degree of Eq(1.2).
\end{lemma}
\begin{proof}
Assume that the point $s_2=(\sigma_2,t_2)\in\Delta$ is an arbitrary point on the root locus of Eq(1.1). The points which satisfy Eq(1.1) are
all on the root locus of Eq(1.1). So, the point $s_2$
surely satisfies Eq(1.1). When $s_2$ is substituted into
Eq(1.1), $K_2G(s_2)\prod^{m}_{l=1}(1-\frac{s_2}{z_{l}})^{\gamma_{l}}G_{lz}(s_2)/\prod^{n}_{j=1}(1-\frac{s_2}{p_{j}})^{\beta_{j}}G_{jp}(s_2)=a_2+ib_2$.
According to Lemma 1.3, here, the gain $K_2$ is a positive real
number, its phase angle of $K_2$ is $2q_1\pi$, the phase angle of
the right side of the equation in this paragraph on the point $s_2$
is $2q_2\pi+\arg{(b_2/a_2)}=2q_2\pi+\alpha_2$.

The factor of the left side of the equation in last paragraph on the
point $s_2$ can be looked as two factors. One is $K_2$, another is
$G(s_2)\prod^{m}_{l=1}(1-\frac{s_2}{z_{l}})^{\gamma_{l}}G_{lz}(s_2)/\prod^{n}_{j=1}(1-\frac{s_2}{p_{j}})^{\beta_{j}}G_{jp}(s_2)$.
The phase angle of the factor of the left side of the equation in
last paragraph on the point $s_2$ is equal to the summation of two
phase angles of two factors of $K_2$  and $G(s_2)\prod^{m}_{l=1}(1-\frac{s_2}{z_{l}})^{\gamma_{l}}G_{lz}(s_2)/\prod^{n}_{j=1}(1-\frac{s_2}{p_{j}})^{\beta_{j}}G_{jp}(s_2)$. The phase angle of the
factor of the left side of the equation in last paragraph on the
point $s_2$ is equal to the phase angle of the right side of the
equation in the last paragraph on the point $s_2$. The phase angle
$2q_2\pi+\alpha_2$ subtracts the phase angle $2q_1\pi$  is equal to
$2q_2\pi+\alpha_2-2q_1\pi=2q\pi+\alpha_2$.

The difference of the phase angle of the factor $a_2+ib_2 $  and the
factor $K_2$ is: $2q\pi+\alpha_2$. So , the phase angle of
expression $G(s_2)\prod^{m}_{l=1}(1-\frac{s_2}{z_{l}})^{\gamma_{l}}G_{lz}(s_2)/\prod^{n}_{j=1}(1-\frac{s_2}{p_{j}})^{\beta_{j}}G_{jp}(s_2)$ is $2q\pi+\alpha_2=2q\pi+\arg{(b_2/a_2)}$. We can obtain:
the phase angle of $W(s)$ is $\varphi(\sigma,t)=\arg(\frac{G_{y}(\sigma,t)}{G_{x}(\sigma,t)})+
\sum^{m}_{l=1}(\gamma_{l}\arg(\frac{t-t_{l}}{\sigma-\sigma_{l}})-\gamma_{l}\arg(\frac{t_{l}}{\sigma_{l}})+
\arg(\frac{G_{zy}(\sigma,t)}{G_{zx}(\sigma,t)}))
-\sum^{n}_{j=1}(\beta_{j}\arg(\frac{t-t_{j}}{\sigma-\sigma_{j}})-\beta_{j}\arg(\frac{t_{j}}{\sigma_{j}})+
\arg(\frac{G_{py}(\sigma,t)}{G_{px}(\sigma,t)}))$. So,
according to the previous proof, we can obtain: $ \varphi(\sigma_2,t_2)=2q\pi+arg{(b_2/a_2)}$,  the point $s_2$ satisfies the
phase angle condition equation (1.2) . Hence, if the point $s_2$ is
on the root locus of Eq(1.1), which its gain is  $K_2$ and its
degree is $2q\pi+\alpha_2$, and satisfies Eq(1.1), then point
$s_2$ satisfies Eq(1.2).
\end{proof}
Lemma 1.5 gives the sufficient condition result of Theorem 1.7. Lemma
1.6 gives the necessary condition result of Theorem 1.7. So, sum up Lemma
1.5 and Lemma 1.6, we can obtain Theorem 1.7.
\begin{theorem}
Let Eq(1.1) be $2q\pi+\alpha$ degree and $s=(\sigma,t)\in\Delta$ be an arbitrary
point. A necessary and sufficient condition
for that the point $s$  is on the root locus of $2q\pi+\alpha$
degree of Eq(1.1) is that point $s$ whether or not satisfies the
phase angle condition equation (1.2).
\end{theorem}

\begin{definition}
Let $\Xi$ be a point set in $\mathbb{C}\cup\{\infty\}$, that are consist of all of points on the path of roots of Eq(1.1) traced out in $\mathbb{C}\cup\{\infty\}$
as $2q\pi+\alpha=2q\pi+\arg(\frac{b}{a})$ is a constant. That path of Eq(1.1) in $\mathbb{C}\cup\{\infty\}$ is called as the root locus of Eq(1.1).
Namely, the set $\Xi$ is the root locus of the $2q\pi+\alpha$ degree.
\end{definition}

\begin{definition}
When the phase angle of meromorphic function $W(s)$ at the left side of the root locus equation (1.1) is $2q\pi+\alpha$ degree, we call the root locus equation (1.1) as the $2q\pi+\alpha$ degree root locus equation.
\end{definition}

In $\mathbb{C}\cup\{\infty\}$, the phase angle of meromorphic function $W(s)$ at the left side of the root locus equation (1.1) is the degree of the root locus of the root locus equation (1.1). So, we can give a definition of the degree of the root locus of Eq(1.1).

\begin{definition}
When the phase angle of meromorphic function $W(s)$ at the left side of the root locus equation (1.1) is $2q\pi+\alpha$ degree, we call the root locus of the root locus equation (1.1) as $2q\pi+\alpha$ degree root locus.
\end{definition}

\begin{lemma}
For all finite zeros of Eq(1.1), they are on the root locus of all degrees of Eq(1.1) from $2q\pi$ degree to $2q\pi+2\gamma_{l}\pi$ degree.
\end{lemma}
\begin{proof}
Assuming that the point $zp$ is an arbitrary finite zeros of Eq(1.1). When the point $zp$ is substituted into meromorphic function $W(s)$ at the left side of Eq(1.1),
we have: $G(zp)\frac{\prod^{m}_{l=1}(1-\frac{zp}{z_{l}})^{\gamma_{l}}G_{lz}(zp)}
{\prod^{n}_{j=1}(1-\frac{zp}{p_{j}})^{\beta_{j}}G_{jp}(zp)}=0$.
This equation can also be expressed as: $G(zp)\frac{\prod^{m}_{l=1}(1-\frac{zp}{z_{l}})^{\gamma_{l}}G_{lz}(zp)}
{\prod^{n}_{j=1}(1-\frac{zp}{p_{j}})^{\beta_{j}}G_{jp}(zp)}=K_{zp}e^{i\theta_{zp}}=0$.
In this equation, $K_{zp}$ is the modulus of the function $W(s)$ and $\theta_{zp}$ is the phase angle of the function $W(s)$.
So, $K_{zp}=0$, and if $\theta_{zp}=2q\pi+\alpha$, $K_{zp}e^{i(2q\pi+\alpha)}=0*(cos(\theta_{zp})+isin(\theta_{zp}))=0$ is true,
and $\theta_{zp}$ is an arbitrary degree from $2q\pi$ degree to $2q\pi+2\gamma_{l}\pi$ degree.

For the point $zp$ that lets the left side of Eq(1.1) obtain 0, no matter what phase angle it is, since its modulus is 0, $cos(\theta_{zp})+isin(\theta_{zp})$ is the non-zero and non-infinity, the modulus of $cos(\theta_{zp})+isin(\theta_{zp})$ is 1. $\theta_{zp}$ represents an arbitrary degree, which shows: no matter what value $\theta_{zp}$ is,
there is $K_{zp}e^{i(2q\pi+\alpha)}=0*(\cos(\theta_{zp})+i\sin(\theta_{zp}))=0$.

This proves: For the finite zero $zp$ of Eq(1.1), when it is substituted into the meromorphic function $W(s)$, the phase angle of the complex value of the meromorphic function $W(s)$ can be the arbitrary $2q\pi+\alpha$ degree. According to Definition 1.10, the phase angle of the meromorphic function $W(s)$ is namely the phase angle of the root locus of Eq(1.1). Because $2q\pi+\alpha$ is an arbitrary degree, when it obtains all of degrees, this shows the finite zeros of Eq(1.1) are simultaneously on the root locus of all of degrees of Eq(1.1) from $2q\pi$ degree to $2q\pi+2\gamma_{l}\pi$ degree.
\end{proof}

\begin{lemma}
For all finite poles of Eq(1.1), they are on the root locus of all degrees of Eq(1.1) from $2q\pi$ degree to $2q\pi+2\beta_{j}\pi$ degree.
\end{lemma}
\begin{proof}
Assuming that the point $pz$ is an arbitrary finite pole of Eq(1.1). When the point $pz$ is substituted into the meromorphic function $W(s)$, we can obtain a value $Gpz=G(pz)\frac{\prod^{m}_{l=1}(1-\frac{pz}{z_{l}})^{\gamma_{l}}G_{lz}(pz)}
{\prod^{n}_{j=1}(1-\frac{pz}{p_{j}})^{\beta_{j}}G_{jp}(pz)}=\infty$.

In complex analysis, 0 can be written as: $0*e^{i\theta}=0*(\cos\theta+i\sin\theta)$.
The points of non-zero and non-infinity finite values can also be written as: $k*e^{i\theta}=k*(\cos\theta+i\sin\theta)$.
The infinity can also be written as: $(+\infty)*e^{i\theta}=(+\infty)*(\cos\theta+i\sin\theta)$. In which, $(\cos\theta+i\sin\theta)$ is a non-zero and non-infinity, and its modulus is 1 of the unit complex number value.
For $(+\infty)*e^{i\theta}=(+\infty)*(\cos\theta+i\sin\theta)$, no matter what value $\theta=2q\pi+\alpha$ obtains on the unit circle in $\mathbb{C}\cup\{\infty\}$, the values of this expression all obtain the infinity.

The above equation can also be expressed as: $Gpz=G(pz)\frac{\prod^{m}_{l=1}(1-\frac{pz}{z_{l}})^{\gamma_{l}}G_{lz}(pz)}
{\prod^{n}_{j=1}(1-\frac{pz}{p_{j}})^{\beta_{j}}G_{jp}(pz)}=K_{pz}e^{i\theta_{pz}}=\infty$.
In this equation, $K_{pz}$ is the modulus of the meromorphic function $W(s)$ and $\theta_{pz}$ is the phase angle of the meromorphic function $W(s)$.
So, $K_{pz}=+\infty$, and if $\theta_{pz}=2q\pi+\alpha$, $K_{pz}*e^{i\theta_{pz}}=(+\infty)*(\cos\theta_{pz}+i\sin\theta_{pz})=\infty$ is true,
and $\theta_{pz}=2q\pi+\alpha$ is an arbitrary degree from $2q\pi$ degree to $2q\pi+2\beta_{j}\pi$ degree.

For the point $pz$ that lets the left side of Eq(1.1) obtain $\infty$, no matter what phase angle it is, since its modulus is $\infty$, $cos(\theta_{pz})+isin(\theta_{pz})$ is the non-zero and non-infinity, the modulus of $cos(\theta_{pz})+isin(\theta_{pz})$ is 1.
$\theta_{pz}=2q\pi+\alpha$ represents an arbitrary degree, which shows: no matter what value $\theta_{pz}=2q\pi+\alpha$ is, there is $K_{pz}e^{i\theta_{pz}}=(+\infty)*(\cos(\theta_{pz})+i\sin(\theta_{pz}))=\infty$.

For the finite pole $pz$ of Eq(1.1), when it is substituted into the meromorphic function $W(s)$, the phase angle of the complex value of the meromorphic function $W(s)$ can be $2q\pi+\alpha$ degree. According to Definition 1.10, the phase angle of the meromorphic function $W(s)$ is namely the phase angle of the root locus of Eq(1.1).
Because $2q\pi+\alpha$ is an arbitrary degree from $2q\pi$ degree to $2q\pi+2\beta_{j}\pi$ degree, when it obtains all of degrees, this shows the finite poles of Eq(1.1) are simultaneously on the root locus of all of degrees of Eq(1.1) from $2q\pi$ degree to $2q\pi+2\beta_{j}\pi$ degree.
\end{proof}

\begin{lemma}
All of the root locus of arbitrary $2q\pi+\alpha$ degree of Eq(1.1) are originated from poles of Eq(1.1),
and are finally received by zeros of Eq(1.1).
\end{lemma}
\begin{proof}
The arbitrary $2q\pi+\alpha$ degree root locus of Eq(1.1) are the curves in $\mathbb{C}\cup\{\infty\}$. So, they all need to have their own origination points and receiving points. According to the relationship between the points in $\mathbb{C}\cup\{\infty\}$ and the root locus of Eq(1.1), the points in $\mathbb{C}\cup\{\infty\}$ can be divided into four types, one is the poles of Eq(1.1), one is the zeros of Eq(1.1), the other is the general finite points on the $2q\pi+\alpha$ degree root locus of Eq(1.1), and the last type is the infinity points in $\mathbb{C}\cup\{\infty\}$.

Except the finite and infinite zeros of Eq(1.1), and except the finite and infinite poles of Eq(1.1), a general finite point in $\mathbb{C}\cup\{\infty\}$ is on a root loci of a certain degree of Eq(1.1). And a general finite point in $\mathbb{C}\cup\{\infty\}$ only has the phase angles of only one degree and non-zero finite gain value, it can not be the originating point and receiving point of the root locus of Eq(1.1).

When the infinite point in $\mathbb{C}\cup\{\infty\}$ let the meromorphic function $W(s)$  obtain a non-zero finite value $A$.
The gain values of the meromorphic function $W(s)$ on the infinite point in $\mathbb{C}\cup\{\infty\}$ are the non-zero finite value $\abs{A}$.
The phase angles of the meromorphic function $W(s)$  are all equal to $2q\pi+arg{(A)}$ degree.
So, these infinite points in $\mathbb{C}\cup\{\infty\}$ are also the ordinary points of the finite values with the gain $\abs{A}$ on the $2q\pi+arg{(A)}$ degree root locus, they can not originate or receive the root locus.

When the infinite points in $\mathbb{C}\cup\{\infty\}$ are the infinite zeros or poles of Eq(1.1). On the infinite zero or pole, the gain value is $K=+\infty$ or $K=0$ respectively.

The finite poles and zeros of Eq(1.1) are on the root locus of all of degrees of Eq(1.1), and the infinite number of the different degrees root locus are at the same one point, so, they satisfy the condition that the root locus can be originated or received. On the finite and infinite poles of Eq(1.1), there is $K=0$. On the finite and infinite zeros of Eq(1.1), there is $K=+\infty$. Thus, we can let the finite and infinite poles of Eq(1.1) as the origination points of the root locus. The finite and infinite zeros of Eq(1.1) are the receiving points of the root locus.
\end{proof}
\begin{definition}
In $\mathbb{C}\cup\{\infty\}$, if an angle $\phi$ is between the positive
direction of the real axis and the tangent line of the root loci of
arbitrary  $2q\pi+\alpha$ degree that is originated from a finite pole,
then angle $\phi$ is called as \emph{angle of origination} of the
root loci of $2q\pi+\alpha$ degree of Eq(1.1).
\end{definition}
\begin{definition}
In $\mathbb{C}\cup\{\infty\}$, if an angle $\varphi$ is between the
positive direction of the real axis and the tangent line of the root
loci of arbitrary  $2q\pi+\alpha$ degree that receives at a finite
zero, then angle $\varphi$ is called as \emph{angle of receiving} of
the root loci of $2q\pi+\alpha$ degree of Eq(1.1).
\end{definition}
In which, when the rotation is counterclockwise, the obtained angle is the
positive angle, when the rotation is clockwise the obtained angle is
the negative angle.
\begin{theorem}
In $\mathbb{C}\cup\{\infty\}$, let $P_k$  be a $\beta_k$ repeated finite poles of Eq(1.1). The angle of origination from $P_k$ on the root locus of arbitrary $2q\pi+\alpha$ degree of Eq(1.1) is:

$\theta_{P_k}=\frac{1}{\beta_k}(2q\pi-\alpha+arg(G(P_k))+\beta_{k}arg(-p_k)-arg(G_{kp}(P_k))+\sum\limits_{l=1}^{m}
(\gamma_{l}arg{(P_k-z_{l})}-\gamma_{l}arg(-z_{l})+arg(G_{lz}(P_k)))-\sum\limits_{j=1,j \neq k}^{n}
(\beta_{j}arg{(P_k-p_{j})}-\beta_{j}arg(-p_j)+arg(G_{jp}(P_k)))$.

In which,  $\beta_k$  is a real number.
\end{theorem}
\begin{proof}
Selecting a point $s_{1}$ on the root locus of $2q\pi+\alpha$
degree, the point $s_{1}$ should be infinitely approaching to $\beta_k$
repeated finite pole $P_k$ that needs to compute its angle of
origination. Because point  $s_{1}$ is infinitely approaching to the
finite pole  $P_k$, the phase angle $\beta_k arg{(s_1-p_j)}$ of vectors between the point $s_{1}$ and poles ($p_j$, $j\neq k$. Namely, except pole  $P_k$.) all can be substituted by angles  $\beta_{k}arg{(P_k-p_{j})}$, $j\neq k$. The phase angle  $\gamma_{k}arg{(s_1-z_{l})}$ of vectors between the point $s_{1}$  and all of the finite zeros $z_l$ of Eq(1.1) all can be substituted by angles $\gamma_{k}arg{(P_k-z_{l})}$.

The point $s_{1}$ is on the root locus of $2q\pi+\alpha$  degree,
because this point is infinitely approaching to the finite pole
$P_k$, the point $s_{1}$ is infinitely approaching to the tangent
line of the root loci of arbitrary $2q\pi+\alpha$  degree that
originates from the finite pole  $P_k$, so, the vector between point
$s_{1}$  and the finite pole $P_k$ is infinitely approaching to the
tangent line of the root loci of arbitrary $2q\pi+\alpha$ degree
that originates from the finite pole $P_k$. So, according to the
Definition 1.14, we can let the phase angle of vectors between point
$s_{1}$ and the finite pole $P_k$ be infinitely approaching to angle
of origination $\theta_{P_k}$.

The point $s_{1}$ must satisfy the phase angle condition equation,
$arg(G(s_1))+\sum\limits_{l=1}^{m}
(\gamma_{l}arg{(s_1-z_{l})}-\gamma_{l}arg(-z_{l})+arg(G_{lz}(s_1)))-\sum\limits_{j=1}^{n}
(\beta_{j}arg{(s_1-p_{j})}-\beta_{j}arg(-p_j)+arg(G_{jp}(s_1)))=2q\pi+\alpha$.
 So, when the finite pole $P_k$
substitutes point  $s_{1}$, we can obtained. $arg(G(P_k))+\sum\limits_{l=1}^{m}
(\gamma_{l}arg{(P_k-z_{l})}-\gamma_{l}arg(-z_{l})+arg(G_{lz}(P_k)))-\sum\limits_{j=1,j \neq k}^{n}
(\beta_{j}arg{(P_k-p_{j})}-\beta_{j}arg(-p_j)+arg(G_{jp}(P_k)))-\beta_k\theta_{P_k}+\beta_karg(-p_k)-arg(G_{kp}(P_k))=2q\pi+\alpha$.

After terms are moved, obtain: $\beta_k\theta_{P_k}=arg(G(P_k))+\sum\limits_{l=1}^{m}
(\gamma_{l}arg{(P_k-z_{l})}-\gamma_{l}arg(-z_{l})+arg(G_{lz}(P_k)))-\sum\limits_{j=1,j \neq k}^{n}
(\beta_{j}arg{(P_k-p_{j})}-\beta_{j}arg(-p_j)+arg(G_{jp}(P_k)))+\beta_karg(-p_k)-arg(G_{kp}(P_k))-2q\pi-\alpha$

The $2q\pi$  and $-2q\pi$ are equivalent, so, $2q\pi$ is substituted
by  $-2q\pi$, after transpositions, we can obtain the formula in Theorem 1.16.
\end{proof}
\begin{theorem}
In $\mathbb{C}\cup\{\infty\}$, let $Z_k$  be a $\gamma_k$ repeated zeros of Eq(1.1). The angle of receiving by $Z_k$ on the root locus of arbitrary $2q\pi+\alpha$ degree of Eq(1.1) is:

$\theta_{Z_k}=\frac{1}{\gamma_k}(2q\pi+\alpha-arg(G(Z_k))+\gamma_{k}arg(-z_k)-arg(G_{kz}(Z_k))-\sum\limits_{l=1, l \neq k}^{m}
(\gamma_{l}arg{(Z_k-z_{l})}-\gamma_{l}arg(-z_{l})+arg(G_{lz}(Z_k)))+\sum\limits_{j=1}^{n}
(\beta_{j}arg{(Z_k-p_{j})}-\beta_{j}arg(-p_j)+arg(G_{jp}(Z_k)))$.

In which,  $\gamma_k$  is a real number.
\end{theorem}
\begin{proof}
Selecting a point $s_{1}$ on the root locus of $2q\pi+\alpha$ degree. The point $s_{1}$ is infinitely approaching to the $\gamma_k$ repeated zero $Z_k$ that needs to compute its angle of receiving. Because the point $s_{1}$ is infinitely approaching to zero $Z_k$. The phase angle $\beta_k arg{(s_1-p_j)}$ of vectors between point $s_{1}$ and all poles $p_j$ of Eq(1.1) all can be substituted by angles $\beta_{k}arg{(Z_k-p_{j})}$. The phase angle $\gamma_{k}arg{(s_1-z_{l})}$ of vectors between point $s_{1}$ and all zeros ($z_{l}$, $l\neq k$, except zero $Z_k$) of Eq(1.1) all can be substituted by angles $\gamma_{k}arg{(Z_k-z_{l})}$.

The point $s_{1}$ is on the root locus of $2q\pi+\alpha$  degree,
because this point is infinitely approaching to the finite zero
$Z_k$, the point $s_{1}$ is infinitely approaching to the tangent
line of the root loci of arbitrary $2q\pi+\alpha$  degree that is
received by the finite zero  $Z_k$, so, the vector between point
$s_{1}$  and the finite zero $Z_k$ is infinitely approaching to the
tangent line of the root loci of arbitrary $2q\pi+\alpha$ degree
that is received by the finite zero $Z_k$. So, according to the
Definition 1.15, we can let the phase angle of vectors between point
$s_{1}$ and the finite zero $Z_k$ be infinitely approaching to angle
of receiving $\theta_{Z_k}$.

The point  $s_{1}$ must satisfy the phase angle condition equation,
$arg(G(s_1))+\sum\limits_{l=1}^{m}
(\gamma_{l}arg{(s_1-z_{l})}-\gamma_{l}arg(-z_{l})+arg(G_{lz}(s_1)))-\sum\limits_{j=1}^{n}
(\beta_{j}arg{(s_1-p_{j})}-\beta_{j}arg(-p_j)+arg(G_{jp}(s_1)))=2q\pi+\alpha$.
 So, when the finite zero $Z_k$
substitutes point  $s_{1}$, we can obtained. $arg(G(Z_k))+\sum\limits_{l=1, l \neq k}^{m}
(\gamma_{l}arg{(Z_k-z_{l})}-\gamma_{l}arg(-z_{l})+arg(G_{lz}(Z_k)))+\gamma_k\theta_{Z_k}-\gamma_karg(-z_k)+arg(G_{kz}(Z_k))-\sum\limits_{j=1}^{n}
(\beta_{j}arg{(Z_k-p_{j})}-\beta_{j}arg(-p_j)+arg(G_{jp}(Z_k)))=2q\pi+\alpha$.

After terms are moved, obtain: $\gamma_k\theta_{Z_k}=-arg(G(Z_k))-\sum\limits_{l=1, l \neq k}^{m}
(\gamma_{l}arg{(Z_k-z_{l})}-\gamma_{l}arg(-z_{l})+arg(G_{lz}(Z_k)))+\sum\limits_{j=1}^{n}
(\beta_{j}arg{(Z_k-p_{j})}-\beta_{j}arg(-p_j)+arg(G_{jp}(Z_k)))+\gamma_karg(-z_k)-arg(G_{kz}(Z_k))+2q\pi+\alpha$.
After transpositions, we can obtain the formula in Theorem 1.17.
\end{proof}

\begin{theorem}
The difference of degrees of the root locus that a pole emits on two opposite directions is  $\alpha_1-\alpha_2=\beta_k\pi$. In which, $\alpha_1$  is the degree that a pole emits the root locus on one direction,  $\alpha_2$  is the degree that the pole emits the root locus on the opposite direction, $\beta_k$  is the exponent of the pole.
\end{theorem}
\begin{proof}
According to Theorem 1.16, two angles of departure at the $\beta_k$  repeated pole $P_k$  on two opposite directions are:

$\theta_{P_k}=\frac{1}{\beta_k}(2q\pi-\alpha_1+arg(G(P_k))+\beta_{k}arg(-p_k)-arg(G_{kp}(P_k))+\sum\limits_{l=1}^{m}
(\gamma_{l}arg{(P_k-z_{l})}-\gamma_{l}arg(-z_{l})+arg(G_{lz}(P_k)))-\sum\limits_{j=1,j \neq k}^{n}
(\beta_{j}arg{(P_k-p_{j})}-\beta_{j}arg(-p_j)+arg(G_{jp}(P_k)))$.

$\theta_{P_k}+\pi=\frac{1}{\beta_k}(2q\pi-\alpha_2+arg(G(P_k))+\beta_{k}arg(-p_k)-arg(G_{kp}(P_k))+\sum\limits_{l=1}^{m}
(\gamma_{l}arg{(P_k-z_{l})}-\gamma_{l}arg(-z_{l})+arg(G_{lz}(P_k)))-\sum\limits_{j=1,j \neq k}^{n}
(\beta_{j}arg{(P_k-p_{j})}-\beta_{j}arg(-p_j)+arg(G_{jp}(P_k)))$.

Two sides of two equations separately subtract, and obtain:  $-\pi=\frac{1}{\beta_k}(-\alpha_1+\alpha_2)$.
 $-\pi\beta_k=-\alpha_1+\alpha_2$. So,  $\alpha_1-\alpha_2=\beta_k\pi$.
\end{proof}

\begin{theorem}
The difference of degrees of the root locus that a zero receives on two opposite directions is   $\alpha_2-\alpha_1=\gamma_k\pi$, In which,  $\alpha_1$ is the degree that a zero receives the root locus on one direction,  $\alpha_2$   is the degree that the zero receives the root locus on the opposite direction, $\gamma_k$  is the exponent of zero.
\end{theorem}
\begin{proof}
According to Theorem 1.17, two angles of arrival at the $\gamma_k$  repeated zero $Z_k$  on two opposite directions are:

$\theta_{Z_k}=\frac{1}{\gamma_k}(2q\pi+\alpha_1-arg(G(Z_k))+\gamma_{k}arg(-z_k)-arg(G_{kz}(Z_k))-\sum\limits_{l=1, l \neq k}^{m}
(\gamma_{l}arg{(Z_k-z_{l})}-\gamma_{l}arg(-z_{l})+arg(G_{lz}(Z_k)))+\sum\limits_{j=1}^{n}
(\beta_{j}arg{(Z_k-p_{j})}-\beta_{j}arg(-p_j)+arg(G_{jp}(Z_k)))$.

$\theta_{Z_k}+\pi=\frac{1}{\gamma_k}(2q\pi+\alpha_2-arg(G(Z_k))+\gamma_{k}arg(-z_k)-arg(G_{kz}(Z_k))-\sum\limits_{l=1, l \neq k}^{m}
(\gamma_{l}arg{(Z_k-z_{l})}-\gamma_{l}arg(-z_{l})+arg(G_{lz}(Z_k)))+\sum\limits_{j=1}^{n}
(\beta_{j}arg{(Z_k-p_{j})}-\beta_{j}arg(-p_j)+arg(G_{jp}(Z_k)))$.

Two sides of two equations separately subtract, and obtain: $-\pi=\frac{1}{\gamma_k}(\alpha_1-\alpha_2)$.

 $-\pi\gamma_k=\alpha_1-\alpha_2$. So,    $\alpha_2-\alpha_1=\gamma_k\pi$.
\end{proof}

\begin{theorem}
For the root locus which are originated from an arbitrary finite pole $P_k$ of Eq(1.1) in $\mathbb{C}\cup\{\infty\}$, if the angle $\theta_{P_k}$ of origination is a independent variable, then, the degree values of the root locus strictly and monotonously increase clockwise.
\end{theorem}
\begin{proof}
If the angle of origination of the $\alpha_1$ degree root locus of Eq(1.1) on the finite pole $P_k$ is: $\theta_{P_k1}=\frac{1}{\beta_k}(2q\pi-\alpha_1+arg(G(P_k))+\beta_{k}arg(-p_k)-arg(G_{kp}(P_k))+\sum\limits_{l=1}^{m}
(\gamma_{l}arg{(P_k-z_{l})}-\gamma_{l}arg(-z_{l})+arg(G_{lz}(P_k)))-\sum\limits_{j=1,j \neq k}^{n}
(\beta_{j}arg{(P_k-p_{j})}-\beta_{j}arg(-p_j)+arg(G_{jp}(P_k)))$. If the angle of origination of the $\alpha_2$ degree root locus of Eq(1.1) on the finite pole  is: $\theta_{P_k2}=\frac{1}{\beta_k}(2q\pi-\alpha_2+arg(G(P_k))+\beta_{k}arg(-p_k)-arg(G_{kp}(P_k))+\sum\limits_{l=1}^{m}
(\gamma_{l}arg{(P_k-z_{l})}-\gamma_{l}arg(-z_{l})+arg(G_{lz}(P_k)))-\sum\limits_{j=1,j \neq k}^{n}
(\beta_{j}arg{(P_k-p_{j})}-\beta_{j}arg(-p_j)+arg(G_{jp}(P_k)))$. Two equations subtract, obtain: $\theta_{P_k1}-\theta_{P_k2}=\frac{1}{\beta_k}(\alpha_2-\alpha_1)$, so, $\alpha_2-\alpha_1=\beta_k(\theta_{P_k1}-\theta_{P_k2})$.

$\beta_k$ is the exponent of pole $P_k$, it is a real number. $\theta_{P_k1}$ and $\theta_{P_k2}$ are separately the origination angle of the root locus of  $\alpha_1$ degree and $\alpha_2$ degree. The angle of origination is defined as a angle which is between the positive direction of the real axis and the tangent line of the root loci that is originated from a finite pole. Beginning at the positive real axis, rotating counterclockwise to the tangent line of the root loci, a angle of origination can be obtained. So, if $\theta_{P_k1}>\theta_{P_k2}$, $\theta_{P_k1}-\theta_{P_k2}>0$. The angle of origination of the root loci of $\alpha_1$ degree is $\theta_{P_k1}$. The angle of origination of the root loci of $\alpha_2$ degree is $\theta_{P_k2}$. $\alpha_2-\alpha_1>0$,  $\alpha_2>\alpha_1$ is true. So, it can prove: for two root locus which are originated on the two different angles of origination from a same pole, the root loci of the larger degree is at the clockwise direction of the root loci of the less degree.

If $\theta_{P_k1}<\theta_{P_k2}$, $\theta_{P_k1}-\theta_{P_k2}<0$.  $\alpha_2-\alpha_1<0$, $\alpha_2<\alpha_1$ is true. So, it can prove: for two root locus which are originated on the two different angles of origination from a same pole, the root loci of the larger degree is at the clockwise direction of the root loci of the less degree.

So, the degree values of the root locus that are originated from this pole $P_k$ strictly and monotonously increases clockwise.
\end{proof}

\begin{theorem}
For the root locus which are received by an arbitrary finite zero $Z_k$ of Eq(1.1) in  $\mathbb{C}\cup\{\infty\}$, if the angle $\theta_{Z_k}$ of receiving is a independent variable, then its degree values strictly and monotonously decrease clockwise.
\end{theorem}

Let the difference of the degree numbers of a pair of the root locus which begin at pole $P_k$   on opposite direction be a positive, according to Theorem 1.18, we can obtain Theorem 1.22. Let the difference of the degree numbers of a pair of the root locus which begin at pole $Z_k$   on opposite direction be a positive, according to Theorem 1.19, we can obtain  Theorem 1.23.

\begin{theorem}
A necessary and sufficient condition for the  $\beta_k$  repeated pole $P_k$  to be simple is that the difference of the degree numbers of any one pair of the root locus which begin at pole $P_k$   on opposite direction  is a positive integer times of  $\pi$.
\end{theorem}

\begin{theorem}
A necessary and sufficient condition for the  $\gamma_k$  repeated zero  $Z_k$ to be simple is that the difference of the degree number of any one pair of the root locus which end at zero $Z_k$  on opposite direction  is a positive integer times of  $\pi$.
\end{theorem}

$\xi(s)=\frac{1}{2}s(s-1)\Gamma(\frac{s}{2})\pi^{-\frac{s}{2}}\zeta(s)$  is a conjugate function, namely, if  $\xi(s)=u(\sigma+it)+iv(\sigma+it)$, then,  $\xi(\overline{s})=u(\sigma+it)-iv(\sigma+it)$.  $\xi(s)=\xi(1-s)$. So,  $\xi(\frac{1}{2}+it)=\xi(\frac{1}{2}-it)$, namely, if  $\xi(\frac{1}{2}+it)=u(\frac{1}{2}+it)+iv(\frac{1}{2}+it)$, then,  $\xi(\frac{1}{2}+it)=\xi(\frac{1}{2}-it)=u(\frac{1}{2}+it)+iv(\frac{1}{2}+it)$. By the property of the conjugate function, we have:  $\xi(\frac{1}{2}-it)=u(\frac{1}{2}+it)-iv(\frac{1}{2}+it)$. So, $v(\frac{1}{2}+it)=0$, on the critical line, xi-function obtains the real number values.

\bibliographystyle{unsrt}

\end{document}